\documentclass [18pt]{article}
\usepackage {graphicx}

\begin{document}
\begin{center}
\textbf{On the Study of Richard Tom Robert Identity}
\end{center}

{Yeong-Shyeong Tsai}

{Department of Applied Mathematics, National Chung Hsing
University, Taichung, Taiwan}

\textbf{Abstract}

The evolution this article is dramatic. In order to estimate the
average speed of mosquitoes, a simple experiment was designed by
Richard (Lu-Hsing Tsai), Tom (Po-Yu Tsai) and Robert (Hung-Ming
Tsai). The result of the experiment was posted in the science
exhibitions Taichung Taiwan 1993. The average speed of mosquitoes
is inferred by the simple relation $m = K\rho cAt$ , where $m$ is
the number of mosquitoes that are caught by the trap, $\rho $ is
the density of the mosquitoes in the experimental device ( a camp
tent), $c$ is the mean speed of the mosquitoes,$A$ is the effect
area of the trap, $t$ is the time interval, and $K$ is a
proportional constant. In order to find the value of $K$. The
identity $E(l) = \alpha V / A$ is discovered, where $V$ is the
volume of a bounded region in the space $R^n$, $A$ is the boundary
area of the region, $\alpha $ is a proportional constant, $E(l)$
is the expectation of the random paths. By a random path in the
region we mean the path which is traced by a random walker who
passes through the region. In this paper, we will show how to get
the data of $\alpha $ generated by computer. Actually, $\alpha = 1
/ K$. Though the rigorous proof is not shown, a sketch proof will
be shown in this paper. The theoretical values of $K$ are $1 / 2$,
$1 / \pi $, $1 / 4$, $2 / (3\pi )$, $^{\ldots }$, for dimension $n
= 1,2,3,4,...,$\\

\textbf{Introduction}

Although random walk is just a simple mathematical model, it can
be used to describe many physical phenomena such as Brownian
motion and diffusion. There is a very important number called
``Boltzmann constant'' in physics. Historically, scientists used
the random walk model to evaluate this constant for the first time
[4]. The Bertrand paradox can be traced back to the Buffon needle
problem established in 1777 [1]. Since then, the research of
Buffon needle problem has been developed into a branch of
mathematics, called geometric probability. The most important
subject in geometric probability is ``random chords'' [1]. It is
very difficult to design an experiment based on its original
definition because an infinite long straight line is required. To
the best of our knowledge, no one has mentioned the relationship
between the random walk model and the random chord problem. In
this article, we use a simple way to show that random walk problem
is essentially a random chords problem.\\

\textbf{A simple experiment}

In 1993, Hung-Ming Tsai [6] had an idea to estimate the mean speed of
mosquitoes. The idea was: mosquitoes were put in a hemispherical tent with a
diameter of 230 cm, and a trap was set up somewhere in the tent. When a
mosquito touched the trap, it would get killed because there we electrified
grids on the trap. And a electric spark made a sound which was sensed by
human ears. After the mosquitoes were put in the tent uniformly, the switch
of the trap was turned on, the number of sparks was counted, and the time
was recorded. The idea was nothing other than an experiment. By this
experiment, he estimated the mean speed of mosquitoes.

What described above is a crude but interesting experiment. Although the
experiment is easy, the mathematics for analyzing the data obtained from the
experiment seems not trivial. Let $t$ be the time interval in which there
are $m$ mosquitoes that are trapped. Let $\rho $ be the density of
mosquitoes in the tent, let $c$ be the mean speed of the mosquitoes and let
$A$ be the effective area of the trap. In order to analyze the experimental
data, we assume

\begin{equation}
\label{eq1}
m = K\rho cAt
\end{equation}

\noindent where $K$ is a proportional constant. From the
experiment, we can easily get data for $m$, $\rho $, $A$, and $t$.
If we assume that the flying of mosquitoes is a three dimensional
random walk, by applying the theory of geometric probability to
compute the value of $K$ , then the mean speed of mosquitoes can
be estimated from equation (\ref{eq1}). We find that it is
nontrivial to find the theoretical values of $K$. Due to the
Bertrand's paradox and the difficulty, we compute the value of $K$
by computer simulation\\

\textbf{The value of K and a new mathematical identity}

Let $E(l)$be the expectation of the length of random paths in the
bounded region. Without loss the generality, we assume the bounded
region is a sphere. In order to simplify the sketch proof of the
identity, we use discrete model. If it is necessary to treat the
continuous case, then the radius of sphere is taken by limit
process. Let $\tau = E(l) / c$ , where $c$ is the mean speed of
mosquitoes. Then $\tau $ is the mean time interval that a mosquito
stays in the sphere. Consider the sequence of time $t_0 $ , $t_1 $
, \ldots , $t_i $, \ldots , where $t_i = i \cdot \tau $ , $i =
0,1,2,...,$. Imagine that at each time $t_i $ , there are exactly
$N$ (new) mosquitoes get into the sphere and $N$ (old) mosquitoes
get out of the sphere. So, at any time, there are exactly $N$
mosquitoes in the sphere. By the definition of density, we have
$\rho = N / V$ , where $V$ is the volume of the sphere. Now, we
consider that the surface of the sphere is the sensor whose area
is $A$. During any time interval of length $\tau $ , there are
exactly $N$ mosquitoes passing through the boundary and getting
into the sphere. So the boundary ( or the sensor ) senses
$N$mosquitoes, during a time interval of $\tau $, that get into
the sphere and hence senses $Nt / \tau $ mosquitoes in any time
interval of $t$. Since $Nt / \tau $is the $m$ in equation
(\ref{eq1}), $\rho = N / V$, and $\tau = E(l) / c$, equation
(\ref{eq1}) becomes

\begin{equation}
\label{eq2}
K = A / (VE(l))
\end{equation}

The mathematical identity is obtained so easily that we can not almost
believe it. Now, we can prove the identity. Since the system

\noindent
is in equilibrium state, the number of mosquitoes in the bounded

\noindent
region is not changed in any time. Therefore, $\Delta m = \Delta N$.

By so called dimension analysis, we have three freedoms to choose

\noindent
the scale in equation (\ref{eq1}) (a) We can take the time scale so that

\noindent
there are $N$ mosquitoes walk into and walk out the bounded region

\noindent
in a unit time, that is, $N = K\rho cA \cdot 1$ . (b) We can choose the
speed of random walker is 1, that is, $c = 1$ in random walk. (c) By the
definition of density, $\rho = n / V$ the quantity $E(l)A / V$ is
dimensionless. We can choose

\noindent
the unit of the length so that $E(l) = 1$ and hence $E(l) = c$. Therefore,

\noindent
equation (\ref{eq1}) becomes $1 = KAE(l) / V$. Then we have, $E(l) = V / (KA)$ We
have complete the proof of identity (\ref{eq2})

By carefully analyzing, we find two properties of random walk
model. One is the property of homogenous, that is, the properties,
such as the locations and the associated probability, of the
points are translation invariance. The other is the properties of
isotropic, that is, the properties, such as the direction of next
step and the associated probability, is rotation invariance.
Therefore, we find the expectation of the length of random paths
is the same as the expectation of the length of random chords. In
order to demonstrate the fact, we use two dimensional model to
show the result. In Figure 1, there is one random path, from point
$A$ to point $B$ then from point $B$ to point $C$. In Figure 2,
there is one random path, from point $D$ to point $B$ then from
point $B$ to point $E$. The events of these two figures are
equally likely. In the Figure 3, there are two random paths, one
is the same as in Figure 1 and the other is the same as in Figure
2. In the Figure 4, there are two random paths, one is the path
from point $A$ to point $B$ then from point $B$ to point $E$. and
the other is the path from point $D$ to point $B$ then from point
$B$ to point $C$. Using these facts, these two cases of events and
their associated probabilities, events in Figure 3 and the events
in Figure 4, must be equivalent to each other. These arguments
support that the expectation of random path (Figure 3) is the same
as that of random chords (Figure 4).\\

 \textbf{The Theoretical
Values of }$K$

Now, we are going to find the expectation of random chord. Without
losing the generality, we solve the two dimensional problems
instead of n-dimensional problem. Though there are so called the
Bertrand paradox in random chords problem, we are able solve
paradox since the results computer simulation will show the
correct answer of this paradox. This verifying work could not be
done in that time.

We consider a circle with radius $r$. The set of all chords is
well defined. Therefore, we can find the expectation (average) of
the lengths of the set of chords. Since the chords are uniformly
and randomly distributed in the circle. we can find some subset of
these chords by defining a equivalent relation. The parallel
relation between chords is a equivalent relation on the set of
that chords. The equivalent relation partitions this set into
disjoint equivalent classes. Clearly, the expectations of these
equivalent classes must be the same one. In the case of finite
set, if the mean of all disjoint subsets is the same one, then the
set has the same as that of these subsets. Therefore, we will
compute the expectation of the random chords in one of these
equivalent classes. In figure 5, the measure of the distribution
is $dy / (2r)$. Therefore, we have expectation of random chords,
$E(C)$,

\begin{center}
$E(C) = \int_{ - r}^r {2xdy} / (2r)$. (3)
\end{center}

\begin{equation}
\label{eq3} E(C) = \pi r^2 / (2r).
\end{equation}

We interpretate equation (\ref{eq3}) as $E(C) = V_n / V_{n - 1} $
since the diameter is the normal section of the circle. This
interpretation can be generalize to higher dimension cases.

From identity (\ref{eq2}), we have

\begin{equation}
\label{eq4} K = V_{n - 1} / A_n ,
\end{equation}

\noindent where $V_{n - 1} $ is the volume of $S^{n - 1}$ the
normal section of $S^n$ and $A_n $ is the boundary area of $S^n$.
Here, we mean that the sphere contains all its interior points. If
it is necessary, we will use the term the sphere surface which is
$S^n$ to be used in usual sense of mathematician. From equation
(\ref{eq4}), we are able to calculate the theoretical values of
all $K$.\\

\textbf{Discovery of some sequence}

The paper, by H. Chalkley, J. Confield and H. Park, stresses how to estimate
the ratio of volume and area. So do their followers'. But our interest is to
find the values of $K$. It can be shown that the values of $K$ depend on the
dimension the space and that the values of $K$ are not dependent on the
shape or the size of the bounded regions. The ratio $V / A$ of $S^n$ is well
known. If we can find the expectation of random chords, then we are able to
get the theoretical values of $K$.

From equation (\ref{eq1}), We can estimate the value of by computer simulation while
the theoretical value of $K$ can be obtained by equation $K = A / (VE(l))$ .
The values of are $K_1 = 1 / 2$, $K_2 = 1 / \pi $ , $K_3 = 1 / 4$, $K_4 = 2
/ (3\pi )$, \ldots $K_7 = 5 / 32$, \ldots , $K_{10} = 128 / (315\pi )$,
\ldots

Our approach is to study the related topics rather than to solve this
problem separately or independently. In order to convince the readers that
the results we obtain are correct and can be verified by computer
simulation, the computer programs have been run on PC, DEC VAX 9000, and IBM
SP2 for more than 5 years. It is computer that generates a phenomenon which
we have never found. From this phenomenon, we try to find some mathematical
model to fit it. It is quite natural to generalize the results to Riemannian
manifold. And then the sequence $1 / 2$, $1 / \pi $, $1 / 4$, $2 / (3\pi
)$,\ldots ,$5 / 32$,\ldots , $128 / (315\pi )$,\ldots , is discovered. So
far, we have found that K$_{d}$'s depend on the dimension of the space or
manifold only.

It takes time to judge how important the identity , $E(l) = A /
KV$ , is, because there are many branches of mathematics such as
geometry, analysis and probability or measure theory, concerning
with this equation . We think that this is a interesting
problem.\\
We have done the computer simulation of random walk on Riemannian
manifold, the surface of the sphere $S^n$, $n \ge 2$. Almost the
same result is obtained. We will organize another paper to show
how to design the computer program for simulating the random walk
on the surface of sphere $S^n$. We shall simply use Riemannian
manifold $S^n$ in usual sense. In $S^2$, we are shocked by a
simple example. Let the boundary region be defined on the north
pole of the earth $S^2$. The region is a circle of which the
center is the north pole and the radius $R$, $R \le \pi r / 2$,
since $R \ge \pi r / 2$, the region becomes a circle on the south
pole. The chord of the circle is defined as the segment of a
geodesic that contains at least one interior point of the circle,
the bounded region. Of course, the two end points must be on the
circle. Let us use the spherical coordinate system, $\rho $,
$\theta $ and $\phi $. Here $\rho $ is constant since we are
studying the problem on $S^2$. For any $\theta $, $0 < \theta \le
\pi / 2$, there is a circle defined on the north pole. The value
of $V$ is $\int_0^\theta {2\pi r\sin \theta } 'rd\theta '$, that
is, $2\pi r^2(1 - \cos \theta )$. The value of $A$ is $2\pi r\sin
\theta $, The value of $K$ is $1 / \pi $. Therefore, the value of
$E(l)$ or $E(C)$ is $\pi r(1 - \cos \theta ) / \sin \theta $, $0 <
\theta \le \pi / 2$. The circle is the north hemisphere of the
earth when $\theta = \pi / 2$. Then the value of $V$is $2\pi r^2$,
the value of $A$ is $2\pi r$ and the value of $E(l)$ or $E(C)$ is
$\pi r$. We find this is a circle of which all chords have the
same length, $\pi r$. Since we were confused by this circle, we
have spent more than one year's time to reinvestigate our
formulation. Finally, we found the result is correct.\\

\textbf{References}

1. H. Solomon, \textit{Geometric Probability }(Society for Industrial and Applied Mathematics,
Philadelphia, Pennsylvania, 1978), pp. 1-32, 127-172 .

2. H. Chalkley, J. Cornfield, H. Park, A method for estimating
volume-surface ratios, Science, 110, pp. 295-297 (1949).

3. C. Kittel AND H. Kroemer, \textit{Thermal Physics} (W.H. Freeman and Company, San Francisco, ed.
2, 1980), pp. 399-414.

4. R.P. Feynman et al., \textit{Lectures on Physics} (California Institute of Technology, 1963), vol. 1,
section 43-3-section 43-5.

5. K. Stowe , \textit{Introduction to Statistical Mechanics and Thermodynamics }(John Wiley, New York, 1984), pp. 324-327.

6. Hung-Ming Tsai et al., Estimation of the Mean Speed of Mosquitoes,
\textit{Science Exhibition, Taichung, Taiwan} (1993)$.$
\newpage

\begin{figure}[htbp]
\centerline{\includegraphics[width=4.17in,height=3.13in]{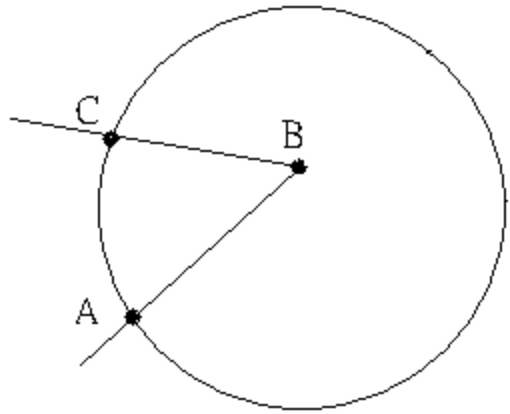}}
\caption{In this figure. There is a mosquito that gets into the
region from boundary point A. Then the mosquito makes a turn at
point B and goes out the region through the boundary point C. }
\label{fig1}
\end{figure}
\newpage

\begin{figure}[htbp]
\centerline{\includegraphics[width=4.17in,height=3.13in]{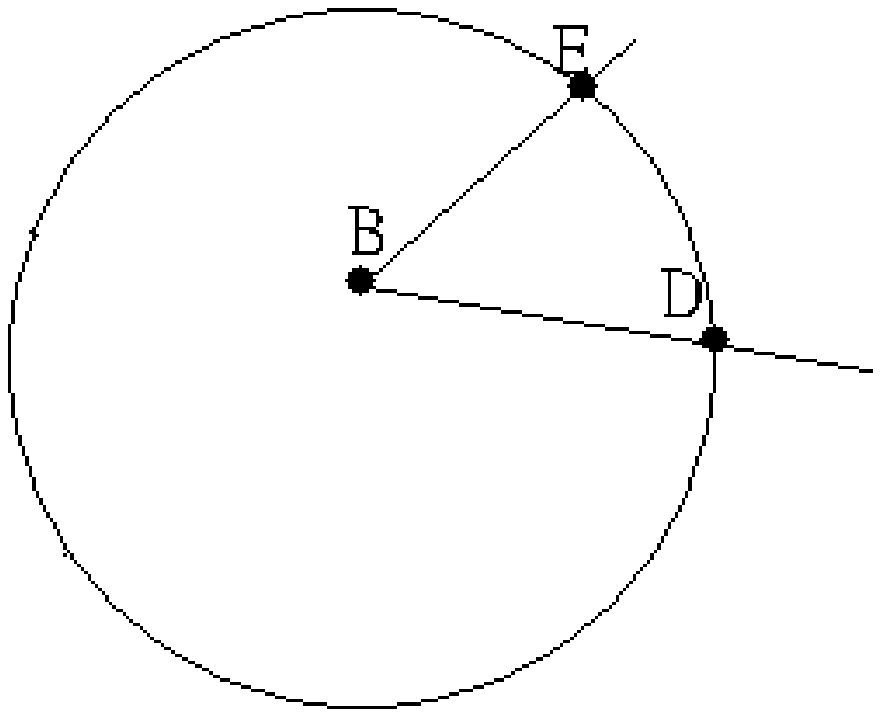}}
\caption{In this figure. There is a mosquito that gets into the
region from boundary point D. Then the mosquito makes a turn at
point B and goes out the region through the boundary point E. }
\label{fig2}
\end{figure}
\newpage

\begin{figure}[htbp]
\centerline{\includegraphics[width=4.17in,height=3.13in]{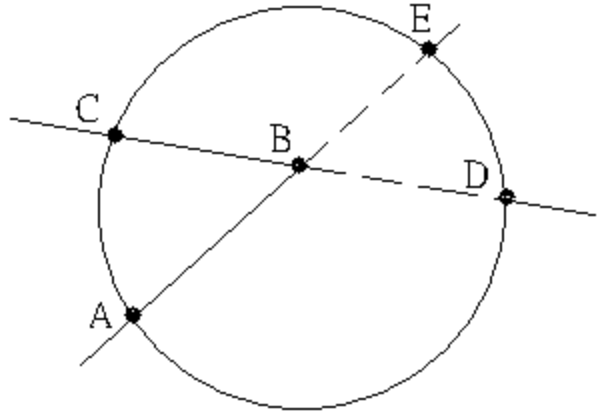}}
\caption{In this figure. There are two mosquitoes that get into
the region. one of them gets into the region from boundary point
A. Then the mosquito makes a turn at point B and goes out the
region through the boundary point C. The other gets into the
region from boundary point D. Then the mosquito makes a turn at
point B and goes out the region through the boundary point E. }
\label{fig3}
\end{figure}
\newpage

\begin{figure}[htbp]
\centerline{\includegraphics[width=4.17in,height=3.13in]{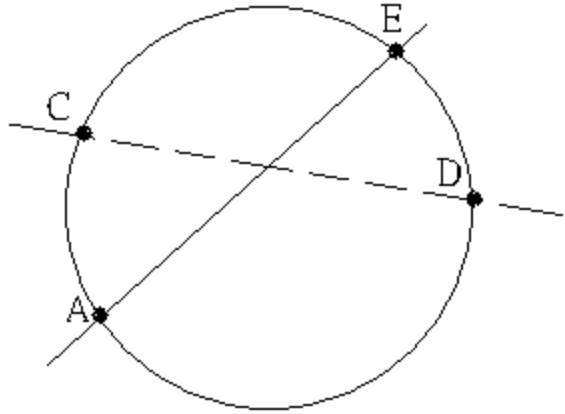}}
\caption{In the figure, there are two random chords, $ AE $ and $
DC $ in the circle. Clearly, The mean the lengths of these two
chords is the same as that of the two random paths in Figure. 3. }
\label{fig4}
\end{figure}
\newpage
\begin{figure}[htbp]
\centerline{\includegraphics[width=4.17in,height=3.13in]{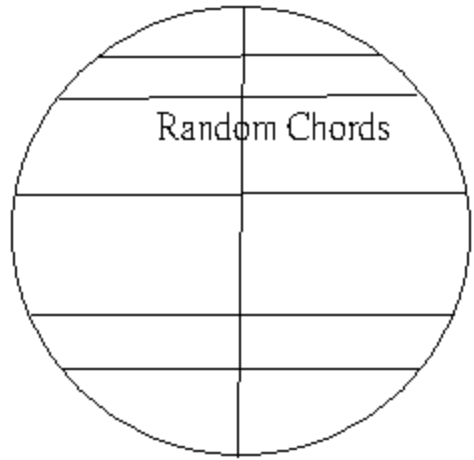}}
\caption{In the figure, some random chords of a class are shown.
The diameter that is perpendicular to the random chords is called
the normal section of the circle. } \label{fig5}
\end{figure}
\newpage

\textbf{the table of data}\

\begin {tabular}{lllll}
\hline
dimension 1 & dimension 2   & dimension 3   & dimension 4   & dimension 5 \\
\hline
0.50065 & 0.313622  & 0.248796  & 0.216937  & 0.188525 \\
0.4994  & 0.319567  & 0.262792  & 0.211288  & 0.18821 \\
0.58505 & 0.313033  & 0.251967  & 0.213581  & 0.197106 \\
0.49935 & 0.319322  & 0.257308  & 0.232106  & 0.185656 \\
0.49895 & 0.323244  & 0.252587  & 0.214341  & 0.18923 \\
0.53095 & 0.316656  & 0.245117  & 0.211194  & 0.19041 \\
0.55985 & 0.346344  & 0.254642  & 0.209172  & 0.195942 \\
0.50035 & 0.3259    & 0.251117  & 0.21765   & 0.189067 \\
0.4995  & 0.314611  & 0.255404  & 0.217041  & 0.193394 \\
0.5001  & 0.314278  & 0.249083  & 0.209794  & 0.188805 \\
0.4988  & 0.324444  & 0.250887  & 0.211703  & 0.193387 \\
0.49885 & 0.321067  & 0.255754  & 0.206331  & 0.18697 \\
0.54215 & 0.3226    & 0.254733  & 0.220672  & 0.190093 \\
0.4994  & 0.320378  & 0.263854  & 0.212503  & 0.195502 \\
0.50035 & 0.313078  & 0.253188  & 0.207934  & 0.188976 \\
0.49935 & 0.321233  & 0.254737  & 0.214747  & 0.196928 \\
0.50005 & 0.3178    & 0.247129  & 0.210894  & 0.190982 \\
0.5003  & 0.3166    & 0.254737  & 0.215356  & 0.186613 \\
0.50065 & 0.349744  & 0.252467  & 0.214447  & 0.195952 \\
0.4988  & 0.314167  & 0.251663  & 0.212528  & 0.18815 \\
1/2 & $1/\pi$  & 1/4   & $2/(3\pi)$   & 3/16 \\
\hline
dimension 6 & dimension 7   & dimension 8   & dimension 9   & dimension 10 \\
\hline
0.171587    & 0.15727   & 0.150828  & 0.141217  & 0.13023 \\
0.170157    & 0.156205  & 0.146288  & 0.140869  & 0.131315 \\
0.170313    & 0.157762  & 0.145272  & 0.136662  & 0.126905 \\
0.172759    & 0.158712  & 0.142739  & 0.137612  & 0.130235 \\
0.174857    & 0.155061  & 0.144548  & 0.14564   & 0.13146 \\
0.178097    & 0.155077  & 0.143702  & 0.140506  & 0.130815 \\
0.170351    & 0.15632   & 0.150119  & 0.134872  & 0.12698 \\
0.169589    & 0.152943  & 0.147578  & 0.137501  & 0.131175 \\
0.173221    & 0.154071  & 0.145909  & 0.136983  & 0.125865 \\
0.176309    & 0.161386  & 0.143639  & 0.13556   & 0.129785 \\
0.170673    & 0.15788   & 0.146711  & 0.133069  & 0.132895 \\
0.167793    & 0.157454  & 0.143225  & 0.136452  & 0.12908 \\
0.169161    & 0.160764  & 0.143214  & 0.140338  & 0.13793 \\
0.170093    & 0.155816  & 0.144177  & 0.141294  & 0.131325 \\
0.175531    & 0.153232  & 0.145216  & 0.135168  & 0.12902 \\
0.169005    & 0.158348  & 0.14707   & 0.13701   &0.12985 \\
0.170316    & 0.153707  & 0.14688   & 0.138528  & 0.13025 \\
0.17744 & 0.158312  & 0.146405  & 0.134812  & 0.138405 \\
0.17088 & 0.160957  & 0.146869  & 0.13817   & 0.131965 \\
0.169431    & 0.157529  & 0.147278  & 0.135202  & 0.1304 \\
$8/(15\pi)$    & 5/32  & $16/(35\pi)$ & 35/256    & $128/(315\pi)$ \\
\hline
\end{tabular}\\

\newpage
\textbf{The computer program}\\

 \noindent
  program sphere

! This main program. id and k are the control variable for

! the dimension of the sace. m is the number

! of partcles. n is the steps of random walk.

\noindent dimension
x(600100,10),dd(600100,10),dx(10),tow(10),tm(10),tn(10)

\noindent open(2,file="D:$\backslash $for$\backslash
$sp41.dat",status="unknown")

\noindent rewind(2)

\noindent write(2,*)"t=0.01 pro=sp n=1000 4th aprali"

\noindent pi=3.1415926

\noindent r=1

\noindent t=0.01

\noindent do id=1,10

\noindent ii=id

\noindent call brn(ii,a,ad)

\noindent nn=a*1.15**id

\noindent n=1000*nn

\noindent m=1000*nn

\noindent do k=1,10

\noindent s=0

\noindent ss=0

! The m number of particles is uniformly distributed.

! The partcles is generated uniforly in the general interval.

! The method of rejction is used to obtain the desired results.

\noindent do i=1,m

100 rs=0

\noindent do j=1,id

\noindent call random{\_}number(xZ)

\noindent x(i,j)=2*xz-1

\noindent rs=x(i,j)**2+rs

\noindent end do

\noindent r1=sqrt(rs)

\noindent if (r1 .GE. r)go to 100

\noindent end do

\noindent do l1=1,n

\noindent ss1=0

! The directions of next step are generated.

\noindent do i=1,m

\noindent call rt(dx,id)

\noindent do j=1,id

\noindent dd(i,j)=dx(j)

\noindent end do

\noindent end do

! The postions of m partcles are updated.

! Of course, the counters s and ss1 are used to

! count the number of partcles that hit the boundary.

! The technique for using two counter s and ss1 insure

! the sum is correct. That is, sss1=ss1+1 is not correc

! counter, when the content is very large.

\noindent do i=1,m

\noindent do j=1,id

\noindent dx(j)=dd(i,j)

\noindent tow(j)=x(i,j)

\noindent x(i,j)=x(i,j)+t*dx(j)

\noindent tn(j)=x(i,j)

\noindent end do

\noindent rs=0

\noindent do j=1,id

\noindent rs=x(i,j)**2+rs

\noindent end do

\noindent r1=sqrt(rs)

\noindent if(r1 .LT. r) go to 300

\noindent call tt(tow,tn,r,dx,id,t)

\noindent do j=1,id

\noindent x(i,j)=tn(j)

\noindent dd(i,j)=dx(j)

\noindent end do

\noindent s=s+1

\noindent ss1=ss1+1

300 end do

\noindent ss=ss+ss1

\noindent ss1=0

\noindent end do

! The theoretical value of K, ty , is computed.

\noindent ii=id

\noindent call voa(ii,a,ad)

\noindent tx=n

\noindent sm=m

\noindent b=a*id

\noindent tk=ss/(sm*b*tx*t)*a

\noindent ty=ad/b

! The computer simulation of K is tk.

\noindent write(2,*)tk,ty

\noindent write(*,*)tk,ty

\noindent end do

\noindent end do

\noindent close(2)

\noindent stop

\noindent end

\noindent subroutine tt(tow,tn,r,dx,id,t)

! In this subprogram, the new postion is determined

! when the partcle collides the boundary of sphere.

! The elastic reflection is computed. A simple property is used.

! That is, any position, including on the boundary, of the
particle is vector which is

! perpendicular to the (tangent of )boundary

\noindent dimension tow(10),tn(10),dx(10),f(10),g(10),h(10)

\noindent ap=0

\noindent pl=0

\noindent do j=1,id

\noindent ap=ap+dx(j)*tow(j)

\noindent pl=pl+tow(j)**2

\noindent end do

\noindent tp=-ap+sqrt(ap**2+1-pl)

\noindent pr=0

\noindent pr1=0

\noindent do j=1,id

\noindent f(j)=tow(j)+tp*dx(j)

\noindent pr=pr+(tn(j)-f(j))*f(j)

\noindent pr1=pr1+dx(j)*f(j)

\noindent end do

\noindent dl=0

\noindent do j=1,id

\noindent tn(j)=tn(j)-2*pr*f(j)

\noindent dx(j)=dx(j)-2*pr1*dx(j)

\noindent dl=dl+dx(j)**2

\noindent end do

\noindent dl=1.0/sqrt(dl)

\noindent do j=1,id

\noindent dx(j)=dx(j)*dl

\noindent end do

\noindent return

\noindent end

\noindent subroutine rt(dx,id)

! In this subprogram a random unit vector is generated for

! direction of next step random walking.

\noindent dimension dx(10)

100 tl=0

\noindent do j=1,id

\noindent call random{\_}number(xz)

\noindent dx(j)=2*xz-1

\noindent tl=tl+dx(j)**2

\noindent end do

\noindent if(tl .EQ. 0) go to 100

\noindent tl=sqrt(tl)

\noindent if(tl .GT. 1) go to 100

\noindent tl=1.0/tl

\noindent do j=1,id

\noindent dx(j)=dx(j)*tl

\noindent end do

\noindent return

\noindent end

\noindent subroutine voa(ii,a,ad)

! In the subprogram, both the volumes of unit sphere and the
normal section

! section of the unit sphere are computed. For example, if input
ii=3,

! then output a is the volume of three dimensional unit sphere and

! the other output ad is the area of a unit two dimensinal sphere,
the disk.

\noindent ii=ii-1

\noindent call brn(ii,a,ad)

\noindent ad=a

\noindent ii=ii+1

\noindent call brn(ii,a,ad)

\noindent return

\noindent end

\noindent subroutine brn(ii,a,ad)

! In this subprogram, the volume of unit sphere is computed

! The parameter ii is the dimension, the parameter a is ouput and

! the parameter ad is a redunduncy

\noindent in=ii/2

\noindent in=ii-2*in+1

\noindent if (in .EQ. 2)go to 400

\noindent f=1

\noindent ls=ii/2

\noindent do j= 1,ls

\noindent f=f*j

\noindent end do

\noindent a=3.1415926**ls/f

\noindent go to 500

400 f=1

\noindent do j=1,ii

\noindent f=f*j

\noindent end do

\noindent l1=(ii-1)/2

\noindent f1=1

\noindent do j=1,l1

\noindent f1=f1*j

\noindent end do

\noindent a=3.1415926**l1*2**ii*f1/f

500 return

\noindent end

\end{document}